\documentclass{amsart}
\RequirePackage{graphicx}

\newtheorem{theorem}{Theorem}[section]
\newtheorem{lemma}[theorem]{Lemma}
\newtheorem{corollary}[theorem]{Corollary}
\newtheorem{proposition}[theorem]{Proposition}
\newtheorem{conjecture}[theorem]{Conjecture}

\theoremstyle{definition}
\newtheorem{definition}[theorem]{Definition}
\newtheorem{example}[theorem]{Example}

\theoremstyle{remark}

\numberwithin{equation}{section}

\newenvironment{case}{\left\{\begin{array}{ll}}{\end{array}\right.}
\newcommand{\fig}[2]{\begin{figure}[tb]\includegraphics*{pjcollins-#1.eps}\caption{#2}\label{fig:#1}\end{figure}}

\newcommand{\exampleref}[1]{Example~\ref{example:#1}}
\newcommand{\figref}[1]{Figure~\ref{fig:#1}}
\newcommand{\figpartref}[2]{Figure~\ref{fig:#1}(#2)}
\newcommand{\secref}[1]{Section~\ref{sec:#1}}
\newcommand{\thmref}[1]{Theorem~\ref{thm:#1}}

\newcommand{\rmand}{\mathrm{\ and\ }}
\newcommand{\rmif}{\mathrm{if\ }}
\newcommand{\calO}{\mathcal{O}}
\newcommand{\calU}{\mathcal{U}}
\newcommand{\ba}{\bar{a}}
\newcommand{\bb}{\bar{b}}
\newcommand{\bc}{\bar{c}}
\newcommand{\R}{\mathbb{R}}
\newcommand{\Q}{\mathbb{Q}}

\newcommand{\N}{\mathbb{N}}

\newcommand{\card}{\#}
\newcommand{\cl}[1]{\mathrm{cl}(#1)}
\newcommand{\fto}{\longrightarrow}
\newcommand{\half}{\frac{1}{2}}
\newcommand{\hgap}{\ \ \ \ \ \ \ \ \ \ }
\newcommand{\homotopic}{\sim}
\newcommand{\htop}{h_\mathit{top}}
\newcommand{\id}{id}

\newcommand{\tendsto}{\longrightarrow}
\newcommand{\us}{{U/S}}

\newcommand{\ind}{\mathrm{Ind}}
\newcommand{\modulo}{\;\mathrm{mod}\;}
\newcommand{\per}{\mathrm{Per}}
\newcommand{\rel}{\mathrm{rel}}
\newcommand{\tr}{\mathrm{Tr}}
\newcommand{\ver}{\mathrm{Ver}}

\newcommand{\rpto}{\fto}
\newcommand{\exto}{\fto}

\newcommand{\reg}{{R}}
\newcommand{\regd}{{\mathbf{R}}}
\newcommand{\regw}{{\mathcal{R}}}
\newcommand{\freg}[1]{{#1_\regd}}

\newcommand{\nieleq}{\simeq}
\newcommand{\ntot}[2]{\bar{N}_{#1}(#2)}
\newcommand{\nsep}[2]{N_{#1}(#2)}
\newcommand{\nyrel}[2]{N^Y_{#1}(#2)}

\newcommand{\cut}{{\mathcal{C}}}
\newcommand{\cutf}[1]{{\cut\!{#1}}}
\newcommand{\cutp}[2]{{\cut_{#1}#2}}
\newcommand{\cuts}[2]{{\cut_{#1}#2}}
\newcommand{\cutt}[1]{{\cut #1}}
\newcommand{\graph}[1]{\mathcal{G}#1}
\newcommand{\perr}[2]{{\per_{#1}(#2)}}
\newcommand{\clperr}[2]{{\widehat{\per}_{#1}(#2)}}

\begin{document}

\title{Dynamics Forced by Surface Trellises}

\author{Pieter Collins}
\address{Department of Mathematics, University of California, Berkeley, California 94704}
\email{piet@math.berkeley.edu}
\thanks{The author wishes to thank Morris Hirsch for his advice and suggestions, which were valuable in writing this paper.} 

\subjclass{Primary 58F15; Secondary 54H25}
\date{\today}

\begin{abstract}
Given a saddle fixed point of a surface diffeomorphism, its stable and unstable curves $W^S$ and $W^U$ often form a homoclinic tangle. 
Given such a tangle, we use topological methods to find periodic points of the diffeomorphism, using only a subset of the tangle with finitely many points of intersection, which we call a {\em trellis}.
We typically obtain exponential growth of periodic orbits, symbolic dynamics and a strictly positive lower bound for topological entropy.
For a simple example occurring in the H\'enon family, we show that the topological entropy is at least $0.527$.
\end{abstract}

\maketitle

%------------------------------------------------------------------------------------------------------------------------------------------------------------
%************************************************************************************************************************************************************
%------------------------------------------------------------------------------------------------------------------------------------------------------------

\section{Introduction}
\label{sec:Introduction}

Let $f:\R^2\fto \R^2$ be a diffeomorphism. 
A fixed point $p$ of $f$ is a \emph{hyperbolic fixed point} if the eigenvalues of $Df(p)$ have modulus $\neq 1$.
By the Stable Manifold Theorem, the stable and unstable sets of $p$ are injectively immersed manifolds, and if $p$ is a saddle point, these manifolds are curves.
If these curves intersect at a point $q$ distinct from $p$, there must be infinitely many intersections, and the stable and unstable curves then form a complicated set called a \emph{homoclinic tangle}.

Homoclinic tangles have been studied extensively, dating back to Poincar\'e and Birkhoff.
 The main result, due in its modern form to Smale, is that a diffeomorphism with a transverse homoclinic point has a horseshoe in some iterate.
While this has been generalised to topologically transverse intersections and quadratic tangencies, little progress has been made in determining more about the actual dynamics forced by a homoclinic tangle.

Since all interesting homoclinic tangles have infinitely many intersection points, we cannot compute them in practice. 
The purpose of this paper is to show that we can obtain interesting information about the dynamics of a system by considering a portion of a homoclinic tangle with only finitely many intersection points. 
We call these objects \emph{trellises}.

We will consider systems on compact surfaces with boundary.
Given a trellis for a system, we find lower bounds for the number of periodic orbits of a given period, and the location of these orbits in terms of the complement of the trellis.
In many cases, we can find a finite type shift which gives a good symbolic description of the system.
The growth rate of the number of periodic points is the same as the entropy of the shift, which is a lower bound for the topological entropy of $f$.
All but finitely many periodic points of the shift are realised by the the original map.

Since all the tools we use are topological, we do not need any differentiability requirements, and we can even weaken the hypothesis that $f$ is invertible.
Further, the methods work equally well for heteroclinic tangles.
We will refer to both homoclinic and heteroclinic tangles as \emph{tangles}.
Note that our terminology differs from that of Easton \cite{Easton86}, who uses the word trellis for what we call a tangle.

Algorithms exist for computing approximations to stable and unstable manifolds for surface diffeomorphisms.
Since transverse intersections of these curves are persistent under perturbations, and trellises contain finitely many intersections, we can often compute trellises precisely.
This allows us to obtain rigorous results about real systems.
As an example, we find symbolic dynamics for the H\'enon map $(x,y)\mapsto(r-x^2+cy,x)$ with parameter values $c=-\frac{4}{5}$ and $r=\frac{3}{2}$, and show that it has topological entropy at least $0.527$.

In \secref{Periodic} we state the definitions and theorems from relative periodic point theory we need to study trellises. 
Proofs and further discussion of the results in this section can be found in \cite{CollinsPPa}.

In \secref{Trellis} we give a formal definition of trellises, and details of the operations we need to study them.
For a trellis $T$, we first cut along the unstable curves of $T$ to obtain a topological pair $\cutt{T}$ consisting of a surface and a subset corresponding to the stable curves.
We then homotopy-retract $\cutt{T}$ onto a graph $\graph{T}$.
If $T$ is a trellis for a map $f$, then we obtain maps $\cutf{f}$ on $\cutt{T}$ and $\graph{f}$ on $\graph{T}$.
We can then use Nielsen theory to show that periodic orbits for the graph map correspond to periodic orbits for the original map $f$.
If the trellis $T$ has transverse intersections and is a subset of a tangle for a homeomorphism $f$, then the growth rate of the number of periodic points of $f$ so found is a lower bound for the topological entropy of $f$.

In \secref{Examples} we give a number of examples showing how we can use these methods to obtain interesting results about the dynamics of maps.

%************************************************************************************************************************************************************

\section{Relative Periodic Point Theory}
\label{sec:Periodic}

In this section we give, without proofs, a brief summary of the definitions and theorems for the relative fixed point theory developed in \cite{CollinsPPa}.
The results are based on standard fixed point theory, a good introduction to which can be found in Brown \cite{Brown71}.

There are two basic types of theory, Lefschetz theory and Nielsen theory.
Both these are homotopy-invariant, and allow for comparison of maps on different spaces.
The Lefschetz theory finds periodic points by looking at cohomology actions on $H^*(X,Y)$, and is most useful when no a priori information about periodic points is available.
The computations involved are similar to those for the cohomological Conley index of Szymczak \cite{Szymczak95FUNDM}, and were motivated by this theory, though some of the topology is complicated since our regions may not have disjoint closures.
The Nielsen theory determines when two periodic points can bifurcate with each other. 
It is most useful when we can explicitly find periodic points for one map in a homotopy class, since we can then decide whether these points exist for other maps in the homotopy class.
When studying trellises, the strongest results are obtained by applying Nielsen theory to maps of divided graphs.

Throughout this section, all topological spaces will be assumed to be compact absolute neighbourhood retracts.
All cohomology groups will be taken over $\Q$.

%------------------------------------------------------------------------------------------------------------------------------------------------------------

\subsection{Topological Pairs, Regions and Itineraries}
\label{sub:toppair}

In this section we define a number of terms which provide a framework for describing dynamics.

\begin{definition}[Topological pairs]
A \emph{topological pair} is a pair $(X,Y)$ where $X$ is a topological space and $Y$ is a closed subset of $X$.
If $(X,Y)$ is a topological pair, we will write $Y^C$ for $X\setminus Y$, the complement of $Y$ in $X$.

A \emph{map of pairs} $f:(X_1,Y_1)\fto(X_2,Y_2)$ is a continuous function $f:X_1\fto X_2$ such that $f(Y_1)\subset Y_2$.
A map of pairs $f:(X_1,Y_1)\fto(X_2,Y_2)$ is \emph{exact} if $f^{-1}(Y_2)\subset Y_1$, or, equivalently, if $f(Y_1^C)\subset Y_2^C$.
\end{definition}

\begin{definition}[Homotopy]
Let $f_0,f_1:(A,B)\fto(X,Y)$.
A \emph{homotopy} from $f_0$ to $f_1$ in the category of topological pairs is a family of maps $f_t:(A,B)\fto(X,Y)$ for $0\le t\le1$ such that such that the function $F:A\times I\fto X$ defined by $F(a,t)=f_t(a)$ is continuous.
We write $f_t:f_0\homotopic f_1$ if $f_0$ is homotopic to $f_1$ via the homotopy $f_t$.
$\homotopic$ induces an equivalence relation on maps of pairs, and we write $[f]$ for the equivalence class of $f$.

A homotopy $f_t$ is a \emph{strong homotopy} if $f_t(a)=f_0(a)$ whenever $f_1(a)=f_0(a)$ an \emph{exact homotopy} if each map $f_t$ is exact.
\end{definition}

\begin{definition}[Regions]
An \emph{region} $\reg$ of a topological pair $(X,Y)$ is an open subset of $X\setminus Y$ such that $\reg\cup Y$ is closed in $X$. 
A \emph{regional space} is a triple $(X,Y;\regd)$ where $(X,Y)$ is a topological pair, and $\regd$ is a set of mutually disjoint regions. 
Note that we do not require $\bigcup\regd$, the union of the regions in $\regd$, to cover $Y^C$.

If $(X_1,Y_1;\regd_1)$ and $(X_2,Y_2;\regd_2)$ are regional spaces, a map $f:(X_1,Y_1;\regd_1)\fto(X_2,Y_2;\regd_2)$ is \emph{region-preserving} if there is a function $\freg{f}:\regd_1\fto\regd_2$ such that for all regions $\reg_1\in\regd_1$, $f(\reg_1)\subset \freg{f}(\reg_1)$, and for all regions $\reg_2\in\regd_2$, $f^{-1}(\reg_2)\subset\bigcup\regd_1$.
\end{definition}

\begin{definition}[Dynamical Systems]
A \emph{dynamical system} on a regional space $(X,Y;\regd)$ is a self-map $f$ of $(X,Y)$.

If $f$ and $g$ are dynamical systems on $(X_1,Y_1;\regd_1)$ and $(X_2,Y_2;\regd_2)$ respectively, a region-preserving map $r:(X_1,Y_1;\regd_1)\fto(X_2,Y_2;\regd_2)$ is a \emph{morphism} from $f$ onto $g$ if there is a map of pairs $s:(X_2,Y_2)\fto(X_1,Y_1)$ such that $r\circ s\homotopic\id$ and $f\homotopic s\circ g\circ r$.
\end{definition}

We interpret $X$ as the base space of the system, $Y$ as invariant set on which the dynamics of $f$ is known, and $\regd$ as the regions in which we are interested in finding symbolic dynamics on.
We will see that if there is a morphism from $f$ onto $g$, then the symbolic dynamics we can compute for $f$ are more complicated than that for $g$.

\begin{definition}[Itineraries and Codes]
\label{defn:itinerary}
Let $f$ be a dynamical system on \linebreak $(X,Y;\regd)$.
A sequence $\reg_0\reg_1\reg_2\ldots$ of regions in $\regd$ is an \emph{itinerary} for $x\in X$ if $f^i(x)\in \reg_i$ for all $i\in\N$.

Let $\per_n(f)$ be the set fixed points of $f^n$ (that is, the set of points of not necessarily least period $n$).
A word $\regw=\reg_0\reg_1\ldots \reg_{n-1}$ on $\regd$ is a \emph{code} for $x\in\per_n(f)$ if $f^i(x)\in \reg_i$ for $0\le i<n$.
We write $\perr{\regw}{f}$ for the set of periodic points with code $\regw$, and $\perr{\regd,n}{f}$ for the set of points with codes in $\regd$ of length $n$.
\end{definition}

Notice that the itinerary is not defined for points which leave $\bigcup\regd$, but since regions are disjoint, it is unique where defined.

%------------------------------------------------------------------------------------------------------------------------------------------------------------

\subsection{Relative Lefschetz Theory}
\label{sec:lefschetz}

Since $X$ and $Y$ are ANRs, we can use the strong excision property to define a \emph{cohomology projection}.

\begin{definition}[Cohomology projection]
Let $\reg$ be a region of $(X,Y)$.
let $j_1:(\reg\cup Y,Y)\hookrightarrow(X,Y)$, $j_2:(X,Y)\hookrightarrow(X,X\setminus \reg)$ and $j_3:(\reg\cup Y,Y)\hookrightarrow(X,X\setminus \reg)$ be inclusions.
$j_3$ is (weakly) excisive, so induces isomorphisms on cohomology. 
The \emph{cohomology projection onto $\reg$} is $\pi_\reg^*=j_2^*\circ(j_3^*)^{-1}\circ j_1^*$.
\end{definition}

Using the cohomology projection, we can restrict the cohomology action of a dynamical system $f$ on $(X,Y;\regd)$ to each region.
Given a word $\regw$ on $\regd$, we can obtain a kind of restricted cohomology action of $f^n$.

\begin{definition}
Let $f$ be a semidynamical system on $(X,Y;\regd)$.
For all $\reg\in\regd$, let $f_\reg^*=\pi_\reg^*\circ f^*$.
For all words $\regw$ on $\regd$ of length $n$, let $f_\regw^*=f_{\reg_0}^*\circ f_{\reg_1}^*\circ\cdots\circ f_{\reg_{n-1}}^*$.
\end{definition}

The Lefschetz number of  $f^*_\regw$ is defined as follows.

\begin{definition}[Lefschetz Number]
\label{defn:lefschetznumber}
The \emph{Lefschetz number} of $f_\regw^*$ is $L(f^*_\regw)=\sum_{i=0}^\infty (-1)^i \tr(f^{(i)}_\regw)$.
\end{definition}

Using this, we can deduce the existence of periodic points with a given code.

\begin{theorem}[Relative Lefschetz Theorem]
\label{thm:rellef}
Let $f$ by a semidynamical system on $(X,Y;\regd)$.
Suppose $\regw$ is a word of length $n$ on $\regd$, and $L(f^*_\regw)\neq 0$.
Then there is a period-$n$ point $x$ such that $x$ is the limit of a sequence $(x_i)$ such that $f^j(x_i)\in \reg_{j\modulo n}$ for all $J<i$.
\end{theorem}

We write $\clperr{\regw}{f}$ for the set of periodic points defined above.
Note that if $x\in\clperr{\regw}{f}$, then, $f^j(x)\in\cl{\reg_{j\modulo n}}$ for all $j$.
We give a result showing how we can compare systems on different spaces.

\begin{theorem}
\label{thm:lefschetzfunctorial}
Let $f$ and $g$ be dynamical systems on $(X_1,Y_1;\regd_1)$ and \linebreak $(X_2,Y_2;\regd_2)$ respectively, and $r$ a morphism from $f$ onto $g$.
Then 
$$\sum_{\regw_1\in\freg{r}^{-1}(\regw_2)}L(f^*_{\regw_1})=L(g^*_{\regw_2})$$
\end{theorem}

%------------------------------------------------------------------------------------------------------------------------------------------------------------

\subsection{Relative Nielsen Theory}
\label{sec:nielsen}

Throughout this section, by \emph{curve} we mean a map $\alpha:(I,J)\fto(X,Y)$, where $I$ is the unit interval.
All homotopies of curves will be relative to endpoints, and we write $\alpha_0\homotopic\alpha_1$ if  $\alpha_0\homotopic\alpha_1$ are homotopic $\rel$ endpoints.

Let $f$ be a dynamical system on $(X,Y;\regd)$, and $n\in\N$.

\begin{definition}
Suppose $x_1,x_2\in\per_n(f)$.
We say $x_1$ is \emph{Nielsen equivalent} to $x_2$, denoted $x_1\nieleq_f x_2$, if there is a subset $J$ of $I$ and exact curves $\alpha_j:(I,J)\exto(X,Y)$ from $f^j(x_1)$ to $f^j(x_2)$ for $j=0\ldots n-1$ such that $\alpha_{j+1\;\modulo\;n}\homotopic f\circ\alpha_j$ for all $j$.
The family $(\alpha_j)$ is a \emph{relating family}.

If $x\in\per_n(f)$, then $x$ is \emph{Nielsen related to $Y$}, denoted $x\nieleq_f Y$ if there is a relating family $(\alpha_j)$ for $x\nieleq_f x$ consisting of exact curves $(I,J)\exto(X,Y)$ for which $J\neq\emptyset$.
If $x\not\nieleq_f Y$, then we say $x$ is \emph{Nielsen separated from $Y$}.
\end{definition}

Clearly $\nieleq_f$ is an equivalence relation.
Equivalence classes of $\per_n(f)$ are called \emph{$n$-Nielsen classes}
We will drop the subscript $f$ where this will cause no confusion.

We have the following important lemma.

\begin{lemma}
\label{lem:yrelated}
If $x_1\nieleq x_2$, then $x_1$ is Nielsen related to $Y$ if and only if $x_2$ is Nielsen related to $Y$.
If $x_1\nieleq x_2$, and $x_1\in\perr{\regw}{f}$, then $x_2\in\perr{\regw}{f}$ or $x_1,x_2\nieleq Y$.
\end{lemma}

We can therefore speak of a Nielsen \emph{class} $Q$ being \emph{Nielsen related to $Y$} or \emph{Nielsen separated from $Y$}.
If $Q$ is Nielsen separated from $Y$, then all points of $Q$ have the same code, which we call the \emph{code for $Q$}.
We let $\nsep{\regw}{f}$ be the number of essential Nielsen classes with code $\regw$, and $\nsep{n}{f}$ the number of Nielsen classes with codes $\regw$ of length $n$.

\begin{theorem}
\label{thm:simopen}
Suppose $Q$ is a Nielsen class of $f$.
Then $Q$ is open in $\per_n(f)$.
\end{theorem}

We can therefore define the index of a Nielsen class $Q$, denoted $\ind(X,Q;f)$ or simply $\ind(Q)$ to be the Lefschetz index $\ind(X,U;f)$, where $U$ is an open neighbourhood of $Q$ containing no other fixed points in its closure.

\begin{definition}[Essential Nielsen class]
A Nielsen class $Q$ is \emph{essential} if \linebreak $\ind(X,Q;f)\neq0$.
\end{definition}

We let $\nsep{n}{f}$, the number of essential Nielsen classes separated from $Y$.
We let $\ntot{n}{f}$ be the total number of essential Nielsen classes, and $\nyrel{n}{f}$ the number of Nielsen classes related to $Y$.
$\nyrel{n}{f}$ may be greater or less than the number of Nielsen classes of $f|_Y$.

The following result is a localisation result for Nielsen theory.

\begin{theorem}
\label{thm:mapchange}
Suppose $f$ and $g$ agree on $\bigcup\regd$.
Then $\nsep{\regw}{f}=\nsep{\regw}{g}$ for all words $\regw$ on $\regd$.
\end{theorem}

If there is a morphism from $f$ to $g$, then $f$ has more Nielsen classes than $g$ in the following sense.

\begin{theorem}
\label{thm:nielsenfunctorial}
Let $f$ and $g$ be dynamical systems on $(X_1,Y_2;\regd_1)$ and \linebreak $(X_2,Y_2;\regd_2)$ respectively, and $r$ morphism from $f$ onto $g$.
Then 
$$\sum_{\regw_1\in\freg{r}^{-1}(\regw_2)}\nsep{\regw_1}{f}\ge\nsep{\regw_2}{g}$$
\end{theorem}

We have the following trivial corollary.

\begin{corollary}
\label{cor:nielhomotopic}
If $g$ is homotopic to $f$, then $\nsep{\regw}{g}=\nsep{\regw}{f}$ for all words $\regw$, and $g$ has at least $\nsep{n}{f}$ points of period $n$.
\end{corollary}

%------------------------------------------------------------------------------------------------------------------------------------------------------------

\subsection{Entropy}
\label{sec:entropy}

There are several ways of defining topological entropy. 
We will use the following definition based on $(\calU,n,f)$-separated sets.

\begin{definition}[Topological entropy]
Let $\calU$ be an open cover of $X$.
Points $x_1,x_2\in X$ are \emph{$(\calU,n,f)$-close} if for all $i<n$ there exist $U_i\in\calU$ such that $f^i(x_1),f^i(x_2)\in U_i$.
Points $x_1,x_2$ \emph{$(\calU,f)$-shadow} each other if they are $(\calU,n,f)$-close for all $n$.

A set $S$ is $(\calU,n,f)$-separated if no two points of $S$ are $(\calU,n,f)$-close.
Let $s(\calU,n,f)$ be the maximum cardinality of a $(\calU,n,f)$ separated set.
Then the \emph{topological entropy} of $f$, written $\htop(f)$ is given by
$$\htop(f)=\sup_\calU\lim_{n\tendsto\infty}\frac{\log s(\calU,n,f)}{n}$$
\end{definition}

We have a classical result that $\htop(f)\ge\limsup_{n\tendsto\infty}\frac{\log N(f^n)}{n}=N_\infty(f)$. 
(See Katok and Hasselblatt \cite{KatokHasselblatt95}).
In other words the growth rate of the number of essential fixed-point classes of $f^n$ is a lower bound for the topological entropy of $f$.

For the relative case, we define the \emph{asymptotic Nielsen number} $N_\infty(f)=\limsup_{n\tendsto\infty}\frac{\log N_n(f)}{n}$.
We would like to show again that $\htop(f)\ge N_\infty(f)$. 
Unfortunately, problems can occur near $Y$, so we introduce an additional hypothesis.

\begin{definition}[Expansive periodicity near $Y$]
\label{defn:expansiveperiodicity}
Let $f$ be a dynamical system on a regional space $(X,Y;\regd)$.
We say $f$ has {\em expansive periodicity near $Y$} if there is a neighbourhood $U_0$ of $Y$ and an open cover $\calU$ of $X$ such that whenever $x_0,x_1\in\perr{\regd,n}{f}\cap W$ are Nielsen separated from $Y$, then either $f^i(x_1)$ and $f^i(x_2)$ are $\calU$-separated for some $i$, or every curve from $x_1$ to $x_2$ in $U_0$ is homotopic to a curve from $x_0$ to $x_1$ which does not intersect $Y$.
\end{definition}

We can show that expansive periodicity near $Y$ is enough to show that the topological entropy is at least the asymptotic Nielsen number.

\begin{theorem}
\label{thm:epte}
Let $f$ be a dynamical system on $(X,Y;\regd)$ with expansive periodicity near $Y$.
Then $\htop(f)\ge \nsep{\infty}{f}$.
\end{theorem}

%************************************************************************************************************************************************************

\section{Trellises}
\label{sec:Trellis}

We now give a formal definition of trellises and two important classes of topological pairs.
We also describe some important operations on these objects.

%------------------------------------------------------------------------------------------------------------------------------------------------------------

\subsection{Trellises}
\label{sec:trellis}

\begin{definition}[Trellis]
A trellis $T$ in a surface with boundary $M$ is a collection $(T^P,T^V,T^U,T^S)$ of subsets of $M\setminus\partial M$ with the following properties.
\begin{enumerate}
\item $T^P$ is finite.
\item $T^U$ and $T^S$ are embedded copies of $T^P\times I$ such that each component of $T^U$ and of $T^S$ contains exactly one point of $T^P$.
\item $T^V=T^U\cap T^S$ is finite.
\end{enumerate}
We write $T=(T^P,T^V,T^U,T^S)$.
\end{definition}

We will write $\us$ for a statement which holds for both the stable ($S$) and unstable ($U$) case.
A trellis is \emph{transverse} if intersections of $T^S$ and $T^U$ are topologically transverse.

\begin{definition}[Segments]
A \emph{segment} is an interval in $T^U$ or $T^S$.
Segments may be open or closed subsets of $T^\us$, or neither.
If $q_1$ and $q_2$ lie in the same component of $T^\us$, we have an \emph{open segment} $T^\us(q_1,q_2)$ and a \emph{closed segment} $T^\us[q_1,q_2]$ between $q_1$ and $q_2$.

An \emph{initial segment} has endpoints $p$ and $q$ where $p\in T^P$.
A \emph{minimal segment} has endpoints $q_1,q_2\in T^V$, and $T^\us(q_1,q_2)$ contains no vertices.
A \emph{maximal segment} has endpoints $q_1,q_2\in T^V$, such that $T^\us[q_1,q_2]$ contains all vertices in that component of $T^\us$.
The \emph{ends} of $T^\us$ are the subsets of $T^\us$ not contained in any maximal segment.
\end{definition}

For our purposes, only the maximal segments of $T^\us$ are important, and so we will sometimes remove the ends of $T^\us$ without explicitly mentioning this.

We now define a natural class of maps between trellises:

\begin{definition}[Trellis Maps]
If $T_1$ is a trellis in $M_1$, $T_2$ is a trellis in $M_2$ and $h:M_1\fto M_2$, we say $h$ is a \emph{trellis map} from $T_1$ to $T_2$ if
\begin{enumerate}
\item $h$ maps $T^P_1$ bijectively with $T^P_2$.
\item $h(T^S_1)\subset T^S_2$.
\item $h^{-1}(T^U_2)\subset T^U_1$.
\end{enumerate}
\end{definition}

Two trellis maps $f_0,f_1$ from $T_1$ to $T_2$ are \emph{homotopic} if there is a homotopy $f_t:f_0\homotopic f_1$ such that each $f_t$ is a trellis map.

The most important trellis maps are those from a trellis $T$ to itself.
If $f:M\fto M$ is such a trellis map, we say \emph{$T$ is a trellis for $f$}.
Clearly, if $f$ is a diffeomorphism with saddle periodic points $T^P$, and stable and unstable curves $T^S$ and $T^U$ with intersection $T^V$, then $(T^P,T^V,T^U,T^S)$ is a trellis for $f$.

We use the more general definition of trellis map to keep a formalism for comparing trellis maps for different trellises; in particular, we have a category of trellises and trellis maps.

%------------------------------------------------------------------------------------------------------------------------------------------------------------

\subsection{Combinatorics of trellises}
\label{sec:combinatorics}

Often the best way of describing a trellis is simply to draw it.
However, it is also useful to have a combinatorial way of describing it. 
We shall only consider the simplest case, namely that of a trellis for a homoclinic tangle on a sphere with transverse intersections.
In this case, $T=(T^P,T^V,T^U,T^S)$, where $T^P$ is a one-point set $\{p\}$, and $T^U$ and $T^S$ are embedded intervals.
We need to choose orientations for $T^U$ and $T^S$.

We now assign coordinates to each point of $T^V$.
The \emph{unstable coordinate} of $q\in T^V$, denoted $n_U(q)$ is $n$ if $q$ is the $n^\mathrm{th}$ point of $T^V$ in the positive direction from $p$ along $T^U$, or the $-n^\mathrm{th}$ point of $T^V$ in the negative direction from $p$.
We define the \emph{stable coordinate} $n_S(q)$ in a similar way.

Merely giving the unstable and stable coordinates of points of $T^V$ is not enough to give a good description of a trellis.
We also need to specify the \emph{orientation} of the crossing of $T^U$ with $T^S$.

The orientation at $q$, written $\calO(q)$ is positive ($+$) if $T^U$ and $T^S$ intersect with the same orientation as they do at $p$, and negative ($-$) if they intersect with the opposite orientation.

We can define a trellis up to ambient isomorphism just by giving $(n_U,n_S,\calO)$ for all points $q\in T^V$.
This description will be called the $(U,S,\calO)$-coordinate description of $T$.

%------------------------------------------------------------------------------------------------------------------------------------------------------------

\subsection{Cutting}
\label{sec:cutting}

Suppose $f:M\fto M$ has trellis $T$.
We would like to obtain a map of pairs from $f$ which captures the action of $f$ on $T$.
The process by which we do this is \emph{cutting} along the unstable curves $T^U$.

\begin{definition}[Cutting]
Let $M$ be a surface.
An embedded curve $\alpha$ is a \emph{cutting curve} if $\alpha\cap\partial M\subset\partial\alpha$.
A finite set of mutually disjoint cutting curves is a \emph{cutting set}.

A surface $\cuts{\alpha}{M}$ is obtained by \emph{cutting $M$ along $\alpha$} if there are curves $\alpha_1,\alpha_2:I\fto\cuts{\alpha}{M}$ in the boundary of $\cuts{\alpha}{M}$ which are disjoint except that we allow $\alpha_1(0)=\alpha_2(0)$ or $\alpha_1(1)=\alpha_2(1)$ (or both), and a map $q_\alpha:\cuts{\alpha}{M}\fto M$ such that $q_\alpha$ is the quotient map for the relation $\alpha_1(t)\sim\alpha_2(t)$, and $\alpha(t)=q_\alpha(\alpha_1(t))=q_\alpha(\alpha_2(t))$.
The quotient map $q_\alpha$ is called the \emph{gluing map}.

If $A$ is a cutting set, we can cut along all curves simultaneously to obtain a surface $\cuts{A}{M}$ and gluing map $q_A$.
\end{definition}

It is a straightforward, though messy, exercise to show that cutting surfaces are unique up to homeomorphism.
Cutting is shown pictorially in \figref{cutting}.

\fig{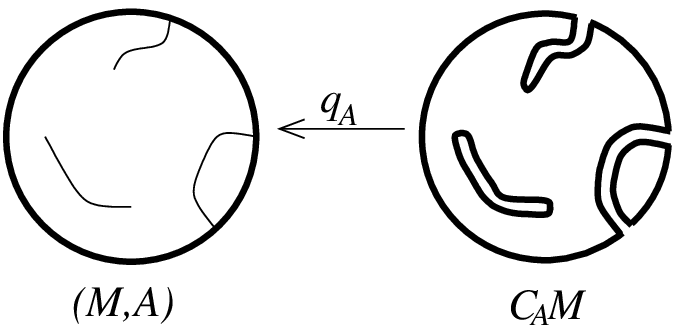}{Cutting along curves}

The gluing map takes $\cuts{A}{M}\setminus q_A^{-1}(A)$ homeomorphically onto $M\setminus A$.
If $x\in A$ then typically $x$ has two preimages under $q_A$, and a neighbourhood $U$ such that $q_A^{-1}(U)$ is homeomorphic to two disjoint copies of the upper-half plane $H$ (and $q_A^{-1}(x)$ lies on the boundaries of these half-planes).
However, if for some arc $\alpha$, $x\in \partial\alpha\setminus\partial M$, then $x$ has a neighbourhood $U$ such that $q_A^{-1}(U)$ is homeomorphic to a single half-plane.

We extend cutting to topological pairs as follows.

\begin{definition}
If $(M,B)$ is a topological pair, and $A$ is a collection of cutting curves, then $\cutp{A}{(M,B)}$ is the pair $(\cuts{A}{M},q_A^{-1}(B))$.
\end{definition}

Given a function $f:M_1\rightarrow M_2$, and cutting sets $A_1$ for $M_1$ and $A_2$ for $M_2$, we would like to know when we can find a map $\cutf{f}:\cuts{A_1}{M_1}\fto\cuts{A_2}{M_2}$ such that $q_{A_2}\circ\cutf{f}=f\circ q_{A_1}$.
The following lemma gives such a condition.

\begin{lemma}
Suppose $M_1$ and $M_2$ are surfaces, $A_1$ and $A_2$ are cutting sets in $M_1$ and $M_2$ respectively, and $f:M_1\fto M_2$ is a map such that $f^{-1}(A_2)\subset A_1$.
Then there is a map $\cutf{f}:\cuts{A_1}{M_1}\fto\cuts{A_2}{M_2}$ such that $q_{A_2}\circ\cutf{f}=f\circ q_{A_1}$.
Further, if $f(B_1)\subset B_2$, then $\cutf{f}(q_{A_1}^{-1}(B_1))\subset q_{A_2}^{-1}(B_2)$
\end{lemma}

\begin{proof}
If $q_{A_1}(x)\in f^{-1}(A_2^C)$, then we can take $\cutf{f}(x)=q_{A_2}^{-1}(f(q_{A_1}(x)))$.
If $f(q_{A_1}(x))$ lies at a point of $A_2$ with one preimages, take $\cutf{f}(x)=q_{A_2}^{-1}(f(q_{A_1}(x)))$.
Otherwise, let $V$ be a neighbourhood of $f(q_{A_1}(x))$ with such that $q_{A_2}^{-1}(V)$ consists of two disjoint copies of $H$.
Let $\hat{U}$ be a semicircular neighbourhood of $x$ such that $q_{A_1}$ maps $\hat{U}$ homeomorphically onto $U$, a subset of $f^{-1}(V)$.
Let $W=U\setminus A_1$.
$W$ is connected, so $f(W)$ is connected, and since $f(W)\subset A_2^C$, $q_{A_2}^{-1}(f(W))$ is connected, so lies in one of the components of $q_{A_2}^{-1}(V)$.
Take $\cutf{f}(x)$ to be the preimage of $f(q_{A}(x))$ under $q_{A_2}$ in this component.

Clearly the map so defined is continuous at $x$, and $\cutf{f}(\cuts{A_1}{B_1})\subset\cuts{A_2}{B_2}$
\end{proof}

Now suppose $T=(T^P,T^V,T^U,T^S)$ is a trellis for a map $f$ on $M$.
We can cut along $T^U$ to obtain a surface $\cuts{T^U}{M}$.
We can also take the preimage of $T^S$ under the gluing map, an obtain a pair $\cutt{T}=(\cuts{T^U}{M},q_{T^U}^{-1}(T^S))$.
For convenience, we will often write $\cutt{T}=(X_T,Y_T)$
An example of the cutting procedure is shown in \figref{cuttrellis}

\fig{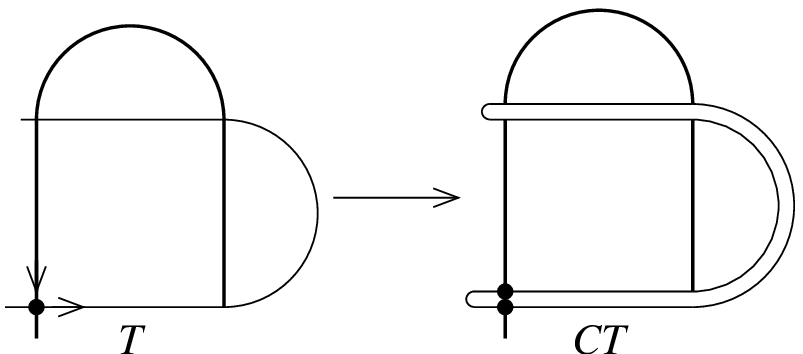}{Cutting along the unstable segment}

Since $f^{-1}(T^U)\subset T^U$, we have a map $\cutf{f}:\cuts{T^U}{M}\fto\cuts{T^U}{M}$, and since $f(T^S)\subset T^S$, $\cutf{f}$ is a map of pairs $\cutf{f}:\cutt{T}\fto\cutt{T}$.
More generally, if $f:M_1\fto M_2$ is a trellis map from $T_1$ to $T_2$, then we can define $\cutf{f}:\cutt{T_1}\fto\cutt{T_2}$.
Since $\cutf{(f\circ g)}=\cutf{f}\circ\cutf{g}$, cutting induces a functor from the trellis category to that of topological pairs.

We now give some trivial, but fundamentally important properties of the $T^U$-cutting projection $q_{A}$.

\begin{proposition}
\label{prop:cuttingproperties}\mbox{\ }
\begin{enumerate}
\item $q_{T^U}$ maps regions of $(M,T^U\cup T^S)$ bijectively with regions of $\cutt{T}$.
\item $f$ has the same periodic orbits as $\cut{f}$, except perhaps for those lying on $T^U$.
\item $q_{T^U}$ is a finite-to-one semiconjugacy, and so $\htop(f)=\htop(\cut{f})$.
\end{enumerate}
\end{proposition}

%------------------------------------------------------------------------------------------------------------------------------------------------------------

\subsection{Cross-Cut Surfaces and Divided Graphs}
\label{sec:collapse}

The relationship between graph maps and surface homeomorphisms has been studied in detail, particularly with regard to Thurston's train tracks and the classification of surface diffeomorphisms. 
More recently, Bestvina and Handel \cite{BestvinaHandell95}, Franks and Misiurewicz \cite{FranksMisiurewicz93} and Los \cite{Los93} produced algorithms for computing the dynamics of isotopy classes of homeomorphisms relative to a finite invariant set.
When studying trellises, we will need to consider \emph{divided graphs}, where we have an invariant subset of the vertex set.
The regions of a divided graph obtained from a trellis are typically very simple (often trees with two or three vertices) making these graphs particularly easy to study.

\begin{definition}[Cross-cut surfaces]
A \emph{cross-cut surface} is a topological pair $(M,A)$, where $M$ is a surface with nonempty boundary, and $A$ is a finite union of disjoint embedded intervals $\alpha$ such that $\alpha\cap\partial M=\partial\alpha$.
$A$ is a \emph{cross-cutting set} and curves $\alpha\in A$ are \emph{cross-cuts}.
\end{definition}

When cutting along $T^U$, all minimal segments of $T^S$ lift to cross-cuts of $\cuts{T^U}{M}$.
If $T$ is a transverse trellis, the endpoints of these lifts are disjoint, so $\cutt{T}$ is a cross-cut surface.

The main property of cross-cut surfaces is that they fibre nicely over graphs.

\begin{definition}[Divided graph]
A divided graph is a topological pair (G,W), where $G$ is a graph (simplicial 1-complex) and $W$ is a subset of $\ver(G)$, the vertex set of $G$.
\end{definition}

We now show that for any pair $(M,A)$ where $M$ is a surface and $A$ consists of nicely embedded curves, there is an exact, homotopy invertible map $r$ to a divided graph.

\begin{theorem}
\label{thm:collapse}
Let $M$ be a surface such that $H_2(M,\emptyset)=0$, and $A\subset M$ set of embedded compact intervals such that $A\cap\partial M$ has only a finite number of components.
Then there is a divided graph $(G,W)$ and an exact map $(M,A)\exto(G,W)$ with a homotopy inverse.
If $M$ is a cross-cut surface, then the homotopy inverse can be made an embedding and all homotopies exact.
\end{theorem}

\begin{proof}
Let $(X,W)$ be the quotient space obtained by collapsing each component of $A$ to a point, and $q$ the quotient map.
Clearly $q$ is exact, and since neighbourhoods of $A$ are topological discs, $q$ has a homotopy inverse $j$.
Further, if $A$ consists of cross-cuts, this homotopy inverse can be made an embedding, as shown in \figref{localretract}

\fig{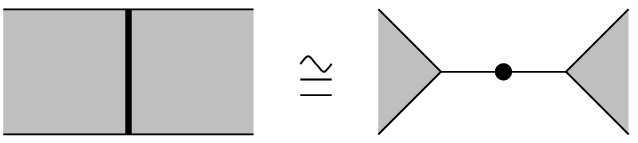}{Exact deformation retract of a cross-cut to a point}

Choose a simplicial subdivision of $X$, such that no simplex contains more than one point of $W$.
Since $X$ is the quotient of a surface by the curves $A$, each 1-simplex of $X$ is contained in no more than two 2-simplexes of $X$.
Then any two vertices lying in the same component of $X\setminus W$ can be joined by an edge-path which does not touch $W$.
Let $Y$ be a minimal 1-complex with the property that any two vertices in the same component of $X\setminus W$ lie in the same component of $Y$.
By the minimality of $Y$, each component of $Y$ is contractible, so $H_2(X,Y\cup W)=0$.
Hence there exists an edge $e$ such that $e\not\in Y$ and $e$ is an edge of exactly one 2-simplex $s$ of $X$.
Let $X_1$ be the simplicial complex formed by removing $e$ and $s$ from $X$.
There is a strong deformation retract $r_1:X\exto X_1$ such that $r_1(s\cup e)\subset \partial s\setminus e$, and both $r_1$ and the corresponding inclusion $i_1$ are exact.
By iterating this procedure to remove one simple at a time, we obtain the graph $(G,W)$.

Since the homotopy inverse for $q$ can be made an exact embedding if $A$ consists of cross cuts, and each inclusion is an exact embedding, we obtain the required homotopy inverse in the case where $A$ consists of cross-cuts.
\end{proof}

Thus there are maps $r:(M,A)\exto(G,W)$ and $s:(G,W)\exto(M,A)$ such that $r\circ s=\id$ and $s\circ r\homotopic\id$.
If $\regd$ is a set of disjoint regions of $(M,A)$, and $\regd_G=\{r(R):R\in\regd\}$, then $r$ is a region-preserving map $(M,A;\regd)\rpto(G,W;\regd_G)$.

Suppose $f$ is a dynamical system on $(M,A;\regd)$.
Let $g=r\circ f\circ s$.
Clearly $r$ is a morphism from $f$ to $g$, so we can study the dynamics of $f$ by studying the dynamics of $g$ using relative Nielsen theory.
If $A$ consists of cross cuts, then since $s\circ g\circ r=s\circ r\circ f\circ r\circ s=s\circ r\circ f\homotopic \id\circ f=f$, there is also a morphism from $f$ to $g$.
In this case, the Nielsen classes of $f$ and $g$ are equivalent.
In the ideal situation, we can find a divided graph $\graph{T}$ and a map $\graph{f}$ such that all periodic points of $\graph{f}$ persist under homotopy.

%------------------------------------------------------------------------------------------------------------------------------------------------------------

\subsection{Graph Maps}
\label{sec:graph}

Under certain conditions, all, or at least all but finitely many, of the periodic points of a system on a graph are unremovable under homotopy.
If there is a morphism from a dynamical system on some other space to such a map, we obtain a lot of information about the periodic points of this system.
One particularly appealing feature of maps on graphs is that we can easily describe homotopy classes combinatorially using simplicial maps.

\begin{definition}
Let $G$ be a graph, $\tilde{G}$ a subdivision of $G$, and $g:\tilde{G}\fto G$ a simplicial map.
We call such a map $g$ a \emph{graph map}.

Let $e$ be an edge of $G$, such that $e=\tilde{e}_1\tilde{e}_2\ldots\tilde{e}_m$, where the $\tilde{e}_i$ are edges of $\tilde{G}$.
Then we write $g(e)=g(\tilde{e}_1)g(\tilde{e}_2)\ldots g(\tilde{e}_m)=e_1e_2\ldots e_n$, the \emph{edge-path action} of $g$.
If $e_{i+1}=\bar{e}_i$ for some $i$, then we say that $g$ \emph{folds} the edge $e$.
\end{definition}

Thus, graph maps either map an edge $e$ to a vertex, or stretch it in a piecewise-linear way over an edge-path $e_1e_2\ldots e_n$ so that the only points of local non-injectivity on $e$ are isolated preimages of vertices.

Dynamics of graph maps can be represented by the \emph{transition matrix}

\begin{definition}[Transition Matrix]
\label{defn:transitionmatrix}
Let $g$ be a graph map of $G$ and let $e_1,\ldots,e_m$ be the edges of $G$.
Let $A$ be the $m\times m$ matrix with $i,j$-th element $a_{ij}$ equal to the number of times $g$ maps edge $e_i$ across $e_j$.
A is the \emph{transition matrix} for $g$.
\end{definition}

If $A$ is the transition matrix for $g$, then we can show that $A^n$ is the transition matrix for $g^n$.
$(A^n)_{ij}$ measures the number of times $g^n$ maps edge $e_i$ across $e_j$.
There must be one periodic point of $g$ of period $n$ in $e_i$ for each time $g^n$ maps $e_i$ across $e_i$ (except in the degenerate case where $g^n(e_i)=e_i$, where all points are periodic by linearity).
Thus there are $(A^n)_{ii}$ period $n$ points of $g$ in $e_i$.

Naively, one would expect $\tr(A^n)=\sum_{i=1}^m(A^n)_{ii}$ to give the total number of points of period $n$ for $g$.
Unfortunately, periodic points in $\ver(G)$ may be counted several times, or not at all.
However, the error between $\tr(A^n)$ and $\card{\per_n(g)}$ is bounded by a constant $c$ independent of $n$.

It is well known that the topological entropy of $g$ is given by the growth rate of the number of periodic points of $g$, $\limsup_{n\rightarrow\infty}\frac{1}{n}\log\tr(A^n)$, and is equal to the Perron-Frobenius eigenvalue of $A$, $\lambda_{\max}(A)$.
$A$ determines a graph with $a_{ij}$ edges from vertex $i$ to vertex $j$, and the dynamics of $g$ are represented by the edge shift on this graph.

Now suppose $(G,W)$ is a divided graph, $\regd$ is a set of disjoint regions, and $g$ is a graph map of $(G,W)$.
We can extend the definition of transition matrices to take into account the regions in $\regd$ as follows:

\begin{definition}
\label{defn:transitionmatrices}
For all regions $R\in\regd$, define an $m\times m$ matrix $P_R$ by $(P_R)_{ii}=1$ if edge $e_i\in R$ and $(P_R)_{ij}=0$ otherwise.
Let $A_R=P_R A$, and $A_\regd=\sum_{R\in\regd}A_R$.
If $\regw$ is a word on $\regd$ of length $n$, let $A_\regw=A_{R_0}A_{R_1}\cdots A_{R_{n-1}}$, the \emph{transition matrix for the code $\regw$}.
\end{definition}

When writing $A_\regd$ we will typically drop rows and columns corresponding to edges not in $\bigcup\regd$, and draw a horizontal line between rows corresponding to edges in different regions.

$\tr(A_\regw)$ is gives the number of points of period $n$ for $g$ with code $\regw$ (except for small errors occurring at vertices).
It is easy to check that
$$\sum_{\regw\in W^m(\regd)}\tr(A_\regw)=\tr(A_\regd^n)\le\tr(A^n)$$

where $W^m(\regd)$ is the set of words on $\regd$ of length $m$.
Again, $\tr(A_\regd^n)$ counts the number of points in $\perr{\regd,n}{g}$, up to an error which is constant in $n$.

We have shown that the periodic points of graph maps are easy to calculate. 
We now define a class of graph maps, called \emph{tight graph maps}, which have minimal dynamics in the homotopy class.

\begin{definition}[Tight Graph Map]
\label{defn:tight}
A graph map $g:(G,W)\fto(G,W)$ is $\regd$-\emph{tight} if for all regions $R\in\regd$, for all edges $e$ in $R$, $g(e)$ does not fold, and if $e_1$ and $e_2$ are distinct edges from the same vertex $v$ in $R\setminus W$, then $g(e_1)$ and $g(e_2)$ have different initial edges.
\end{definition}
 
Not every map of a divided graph is homotopic to a tight graph map, but all the maps of cross-cut surfaces we study are exactly homotopy retract onto a tight graph map, and we conjecture that this is true in general.

The fundamental theorem on tight graph maps is that the periodic points lie in different Nielsen classes, and that, typically, these Nielsen classes are essential.

\begin{theorem}
\label{thm:rexpand}
Suppose $g$ is $\regd$-tight and $x_1,x_2\in\perr{\regd,n}{f}$.
Then either $x_1$ and $x_2$ lie in different Nielsen classes, or there is an edge-path joining $x_1$ to $x_2$ which is fixed by $g^n$.
Further, if $x\in\perr{\regd,n}{f}$, then either $\ind(x;f)\neq0$ or $x\in\ver(G)$. 
\end{theorem}

\begin{proof}
Suppose $x_1\neq x_2$ are Nielsen-equivalent, and $\alpha_j:(I,J)\exto(G,W)$ is a relating family for $x_1\nieleq x_2$.

Suppose $J\neq\emptyset$.
Let $s=\inf J$ and $y=\alpha(s)$
Since $x_1\not\in Y$, $s>0$.
Let $\beta_j:(I,\{1\})\exto(X,Y)$ be given by $\beta_j(t)=\alpha_j(t/s)$.
Then $(\beta_j)$ is a relating family for $x_1\nieleq y$, and further, there are regions $R_j\in\regd$ such that $\beta_j(I)\subset R_j$.

If $J=\emptyset$, then we let $\beta_j=\alpha_j$, so again there are regions $R_j\in\regd$ such that $\beta_j(I)\subset R_j$.

By homotoping if necessary to remove any folds, we can assume that all curves $\beta_j$ are locally injective.
Since $g$ is $\regd$-tight, $g(\beta)$ is locally injective, so, up to parameterisation, $g\circ\beta_j=\beta_{j+1}$
Hence $g^n(\beta_0(I))=\beta_0(I)$ so $g^n\circ\beta_0\homotopic\beta_0$.
Thus $g^n\circ\beta=\beta$ , and so all points of $\beta$ are fixed by $g^n$.

If $x$ is an isolated repelling fixed point of a graph map $f$ and $x$ does not lie on a vertex of $G$, then $\ind(G,x;g^n)=\pm1$.
\end{proof}

%------------------------------------------------------------------------------------------------------------------------------------------------------------

\subsection{Entropy of Trellis Maps}
\label{sec:trellisentropy}

We now show that we can find a lower bound for the entropy of a trellis map in terms of the asymptotic Nielsen number.
By \thmref{epte}, we need only show that $\cutf{f}$ has expansive periodicity near $Y_T$.

\begin{theorem}
\label{thm:trellisentropy}
If $f$ is a homeomorphism with trellis $T$ such that $T^P$ consists of hyperbolic periodic points, then $\cutf{f}:(X_T,Y_T)\fto(X_T,Y_T)$ has expansive periodicity near $Y_T$.
\end{theorem}

\begin{proof}
Since $Y_T$ is the inverse image under the glueing map of a submanifold of stable manifold for $f$, $Y_T$ has a neighbourhood $W$ for which every point of $W\setminus Y_T$ eventually leaves $W$.
Since $Y_T$ is a union of disjoint copies of an interval with endpoints in $\partial X_T$, we can find neighbourhoods $V_1$, $V_2$ and $V_3$ of $Y_T$ each of which deformation retract onto $Y_T$ such that $\cl{V_1}\subset V_2$, $\cl{V_2}\subset V_3$, $\cutf{f}(V_1)\subset V_2$, and every point of $V_1\setminus Y_T$ eventually leaves $V_1$
Choose an open cover $\calU$ containing the components of $V_1$ and $V_2\setminus Y_T$, and such that for all other $U\in\calU$, $U\cap V_1=\emptyset$ and $U$ intersects at most one component of $V_2\setminus Y_T$
(This is where we need $\cl{V_2}\subset V_3$).
Let $U_0=V_1$.
We claim that $\calU$ and $U_0$ are the required open cover and neighbourhood of $Y_T$.

First notice that if $x_1$ and $x_2$ lie in the same component of $V_1$, but different components of $V_1\setminus Y_T$ (equivalently, every path from $x_1$ to $x_2$ in $V_1$ crosses $Y$), then $f(x_1)$ and $f(x_2)$ lie in different components of $V_2$.
Suppose $x_1,x_2\in U_0\setminus Y_T$, and $f^j(x_1)$ and $f^j(x_2)$ are $\calU$-close for all $j$.
Then there exists least $i$ such that either $f^i(x_1)$ or $f^i(x_2)$ are not in $U_0=V_1$.
By minimality of $i$, $f^i(x_1),f^i(x_2)\in V_2$.
Since $f^i(x_1)$ and $f^i(x_2)$ are $\calU$-close, they must lie in the same component of $V_2\setminus Y_T$.
This means that $x_1$ and $x_2$ lie in the same component of $V_1\setminus Y_T$, and since components of $V_1$ are simply connected, every path in $V_1$ from $x_1$ to $x_2$ is homotopic to one which does not intersect $Y_T$.
\end{proof}

We can use this to show that the entropy of a map with trellis $T$ is at least the asymptotic Nielsen number of $\cutf{f}$.

\begin{corollary}
\label{cor:trellisentropy}
If $f$ is a homeomorphism with transverse trellis $T$ such that $T^P$ consists of hyperbolic periodic points, then $\htop(f)\ge N_\infty(\cutf{f})$.
\end{corollary}

\begin{proof}
$\htop(f)=\htop(\cutf{f})$ since the gluing map is a finite-to-one surjective semiconjugacy, and $\htop(\cutf{f})\ge N_\infty(\cutf{f})$ by \thmref{trellisentropy} and \thmref{epte}.
\end{proof}

If the homeomorphism $f$ for the trellis $T$ is clear, we will sometimes call $N_\infty(\cut{f})$ the \emph{entropy of $T$}.

%************************************************************************************************************************************************************

\section{Examples}
\label{sec:Examples}

%------------------------------------------------------------------------------------------------------------------------------------------------------------

\begin{example}[The Smale horseshoe]

First we give a familiar example, the Smale horseshoe map.
Recall that the Smale horseshoe map $f:S^2\fto S^2$ maps the stadium-shaped area of \figref{smale} into itself as shown, mapping the square $S$ linearly across itself with uniform expansion in the horizontal direction and contraction in the vertical direction.
\fig{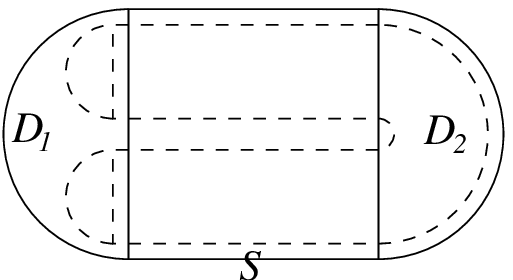}{Smale horseshoe Map}
$f$ maps the semicircular region $D_1$ into itself so that all points in $D_1$ are attracted to a fixed point, and maps $D_2$ into $D_1$.
Outside the stadium, $f$ has a single repelling fixed point.

There is a hyperbolic saddle point in $S$, and the stable and unstable curves form a homoclinic tangle. 
The \emph{horseshoe trellis} $T_2$ is the subset of the tangle in shown in \figpartref{xsmale}{a}.
Except for two fixed points outside $S$, the nonwandering set $\Lambda$ of $f$ lies in the regions $R_1$ and $R_2$.
\fig{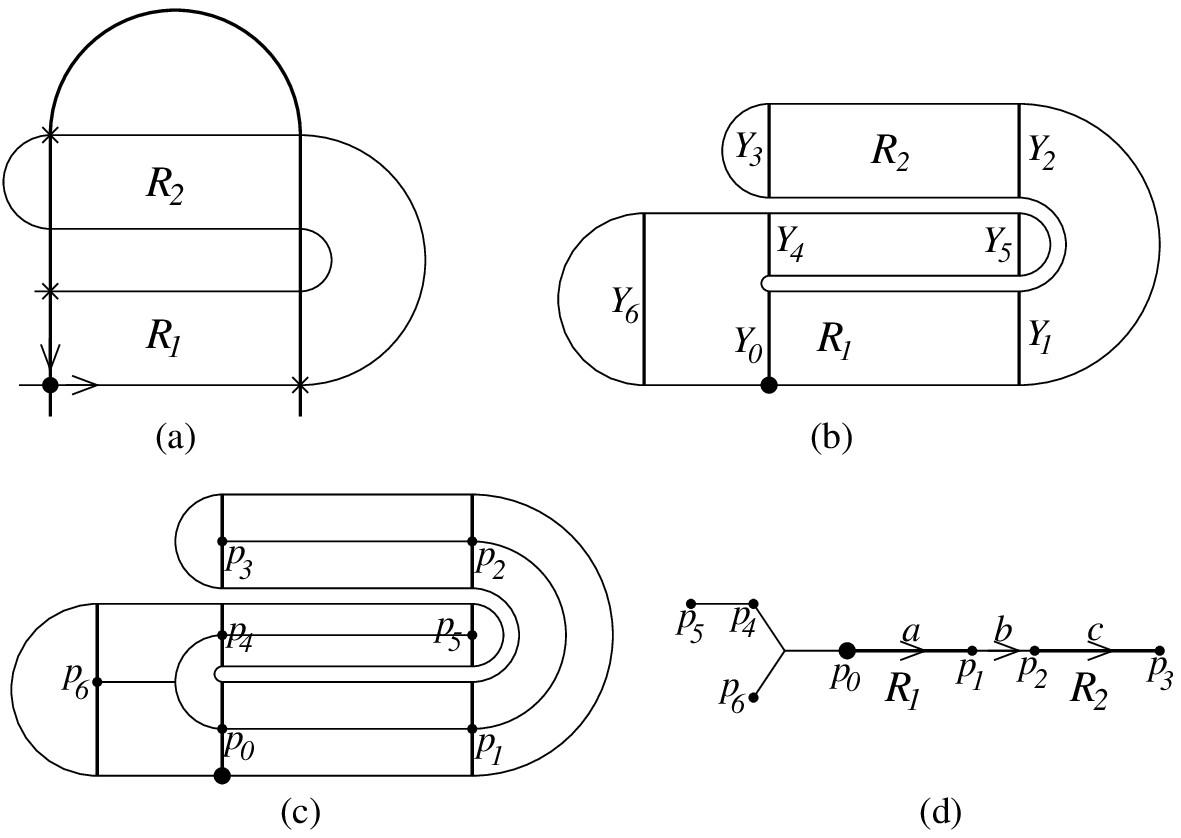}{(a)~Horseshoe trellis $T_2$, (b)~Cut surface $\cutt{T_2}$, (c)~Embedded graph $\graph{T_2}\subset\cutt{T_2}$, (d)~Graph $\graph{T_2}$ with edges labeled}

The $(U,S,\calO)$-coordinates for the vertices are
$$(0,0,+),(1,7,-),(2,4,+),(3,3,-),(4,2,+),(5,5,-),(6,6,+),(7,1,-)$$

To study the dynamics, we first cut along the unstable set $T_2^U$ of the trellis (dropping the ends) as shown in \figpartref{xsmale}{b}.
This gives us a topological pair $\cutt{T_2}=(X_{T_2},Y_{T_2})$, where $X$ is the surface obtained by the cutting, $Y_{T_2}$ is a subset of $X_{T_2}$ corresponding to the stable set $T_2^S$ of the trellis.
$f$ naturally induces a map $\cut{f}$ of $\cutt{T_2}$.

Let $G_{T_2}$ be graph embedded in $\cutt{T_2}$ as shown in \figpartref{xsmale}{c}.
Letting $W_{T_2}=G_{T_2}\cap Y_{T_2}$, we obtain a topological pair $\graph{T_2}=(G_{T_2},W_{T_2})$ onto which we can deformation retract $(X_{T_2},Y_{T_2})$.
This collapsing induces a map $\graph{f}$ on $\graph{T_2}$.

Just by knowing the action of $f$ on $T_2^S$, we can deduce the action of $\graph{f}$ on $W$.
In this case we have 
$$p_0,p_3,p_4\mapsto p_0,\ p_1,p_2,p_5\mapsto p_3 \rmand p_6\mapsto p_4$$
Since $\graph{T_2}$ is a tree, this determines the homotopy class of $\graph{f}$ as a self-map of $\graph{T_2}$ completely.

A tight graph map in the homotopy class of $\graph{f}$ maps the arcs corresponding to regions $R_1$ and $R_2$ across each other.
Using the labeling of \figpartref{xsmale}{d}, we have
$$a\mapsto abc \rmand c\mapsto\bar{c}\bar{b}\bar{a}$$
Thus $\graph{T_2}$ must have a subset on which $\graph{f}$ is conjugate to the one-sided shift on two symbols.
Therefore, the trellis forces dynamics conjugate with the shift on two symbols.
In particular, any map with the same trellis as the Smale horseshoe $f$ must have entropy $\htop\ge\log2$.
\end{example}

%------------------------------------------------------------------------------------------------------------------------------------------------------------

\begin{example}[Iterates of trellis maps]
\label{example:iterate}

Again consider the trellis $T_2$ of and let $f$ be the second iterate of the horseshoe map.
One might expect the homotopy class of $f$ to have \emph{more} entropy than that of $f$.
However, $\graph{f}$ maps all points $p_0\ldots p_6$ to $p_0$ so is homotopic to a constant map.
Thus we obtain no information about the dynamics.
We can find diffeomorphisms homotopic to $f$ with this trellis and arbitrarily small entropy.

\end{example}

%------------------------------------------------------------------------------------------------------------------------------------------------------------

\begin{example}[Trivial trellises]
\label{example:trivial}

Consider the trellis  $T_1$ of \figpartref{xtype1}{a}which is a subset of the horseshoe trellis, and let $f$ be the horseshoe map.
Cutting along the unstable manifolds we obtain the surface $\cutt{T_1}$ shown in \figpartref{xtype1}{b}.
\fig{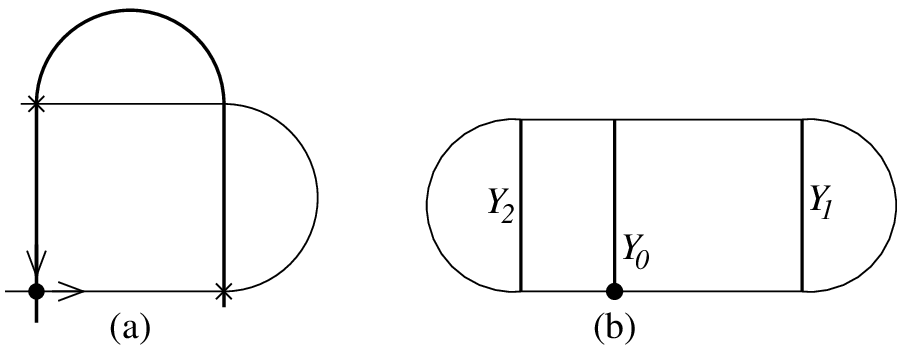}{(a)~Trellis $T_1$, (b)~Surface $\cutt{T_1}$}
The components $Y_0$, $Y_1$ and $Y_2$ of $Y_{T_1}$ all map to $Y_0$ under $\cut{f}$, so $\cut{f}$ is homotopic to a constant.
Therefore, our topological methods give no interesting dynamics.

An even more extreme example is given by the trellis $T_0$ of \figpartref{xtype0}{a}
Cutting along the unstable manifolds we obtain the surface $\cutt{T_0}$ of \figpartref{xtype0}{b}.
\fig{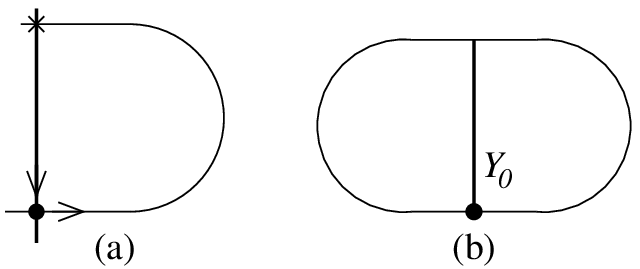}{Trellis $T_0$ and cut surface $\cutt{T_0}$}
All maps on $\cutt{T_0}$ are homotopic to a constant, so again, our topological methods to any map with this trellis yields no information.

In each of these cases, we know that if $f$ is a diffeomorphism with this trellis, $\htop(f)>0$.
However, we can find diffeomorphisms with arbitrarily small entropy.
\end{example}

%------------------------------------------------------------------------------------------------------------------------------------------------------------

\begin{example}[The type-3 trellis]

The type-$3$ trellis $T_3$ is the simplest nontrivial trellis other than the horseshoe.
It trellis occurs in the H\'enon map for a range of parameter values, and a particular case is shown in as in \figref{henon}.
\fig{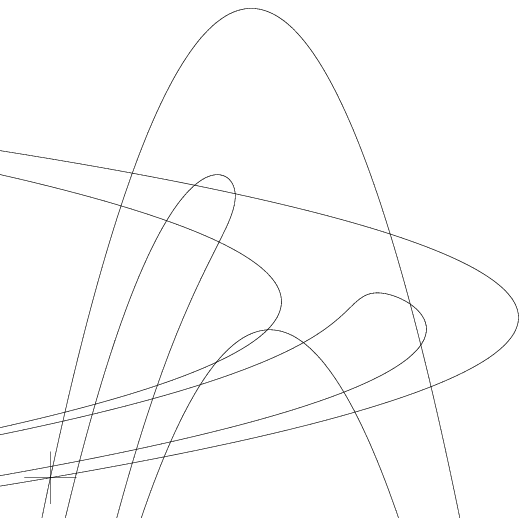}{Trellis in the H\'enon map with parameters $c=-\frac{4}{5}$ and $r=\frac{3}{2}$}
This figure was drawn using the DsTool implementation of the algorithm of Krauskopf and Osinga \cite{KrauskopfOsinga98}.
The trellis is shown in \figpartref{xtype3}{a}.
\fig{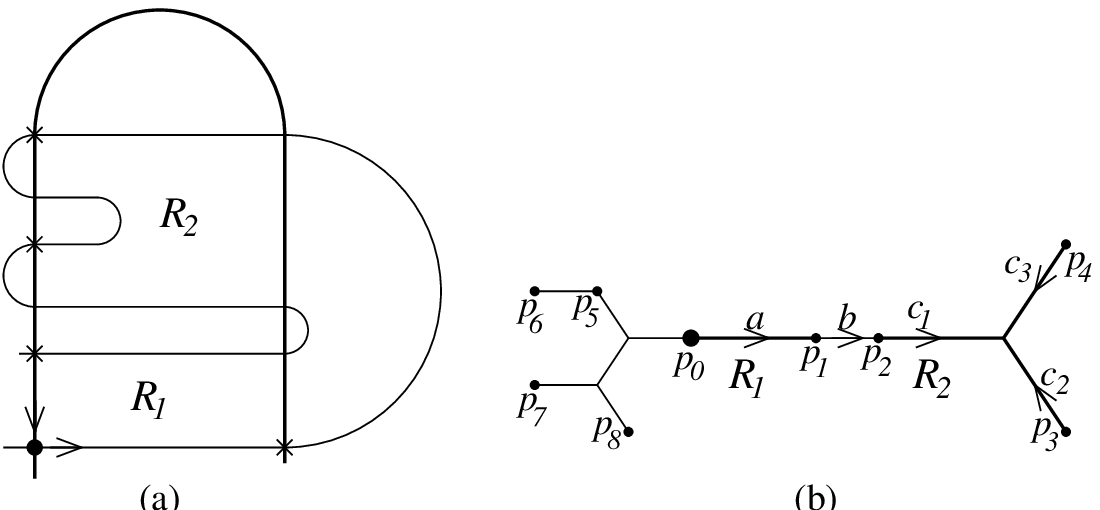}{(a)~Type $3$ trellis $T_3$, (b)~Graph $\graph{T_3}$}
The $(U,S,\calO)$-coordinates for the vertices are
\begin{eqnarray*}
(0,0,+),(1,9,-),(2,6,+),(3,5,-),(4,4,+),(5,3,-),\\
(6,2,+),(7,7,-),(8,8,+),(9,1,-)
\end{eqnarray*}
and the vertices map $(1,9,-)\mapsto(3,5,-)\mapsto(5,3,-)\mapsto(9,1,-)$.

Cutting along the unstable manifold, we obtain the surface $\cutt{T_3}$ and the embedded graph $\graph{T_3}$ as shown in \figpartref{xtype3}{b}.
The action on the distinguished vertex set is 
$$p_0,p_4,p_5\mapsto p_0,\;p_1,p_2,p_6\mapsto p_3,\;p_3\mapsto p_4,\;p_7\mapsto p_5 \rmand p_8\mapsto p_7$$
The graph is a tree, and the regions $R_1$ and $R_2$ are expanding under the the tight map
$$a\mapsto abc_1\bar{c}_2, \; b\mapsto\cdot, \; c_1\mapsto c_2, \; c_2\mapsto c_3 \rmand  c_3\mapsto abc_1 $$
This gives transition matrix (on $\{a,c_1,c_2,c_3\}$)
$$A=\left( \begin{array}{cccc}
      1&1&1&0\\ \hline
      0&0&1&0\\  
      0&0&0&1\\
      1&1&0&0\\
 \end{array} \right)$$
The horizontal line in the matrix separates the rows corresponding to edges of $\reg_1$ from edges of $R_2$.
The edge shift for this transition matrix is given in \figref{shift3}, and since $a\subset R_1$ and $c_1,c_2,c_3\subset R_2$ we obtain a sofic shift on regions.
\fig{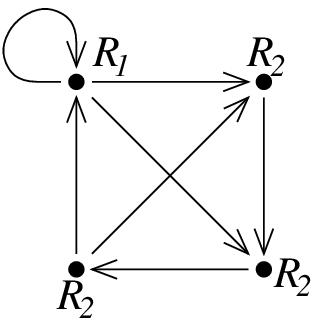}{Shift for $T_3$}
The characteristic polynomial of $A$ is $\lambda(\lambda^3-\lambda^2-2)$, and the Perron-Frobenius eigenvalue $\lambda_{\max}$ of $A$ therefore satisfies $\lambda_{\max}^3-\lambda_{\max}^2-2=0$.
The value of $\lambda_{\max}$ is approximately $1.70$, giving a lower bound of $0.527$ for the topological entropy.
\end{example}

%------------------------------------------------------------------------------------------------------------------------------------------------------------

\begin{example}[The type-$n$ trellis]

The horseshoe trellis and type-$3$ trellis are part of a family of {\em simple trellises}.
The general type-$n$ trellis has vertices with coordinates
\begin{eqnarray*}
(0,0,+),(1,2n+3,-),(2,2n,+),(3,2n-1,-),(4,2n-2,+)\ldots(2n-1,3,-)\\
(2n,2,+),(2n+1,2n+1,-),(2n+2,2n+2,+),(2n+3,1,-)
\end{eqnarray*}

We consider trellis maps taking $(1,2n+3,-)$ to $(3,2n-1,-)$.
The graph $\graph{T_n}$, shown in \figref{xtypen}, has two expanding regions $R_1$ and $R_2$ under the tight map.
\fig{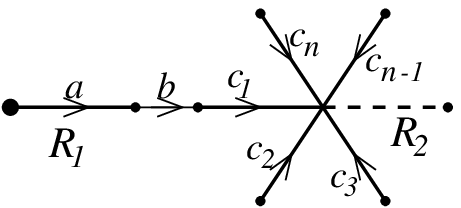}{Graph for $T_n$}
$R_1$ has a single edge $a$, and $R_2$ has edges $c_1,c_2\ldots c_n$ which map:
\begin{eqnarray*}
   a & \mapsto & abc_1\bar{c}_2\\
   c_i & \mapsto & \begin{case} c_{i+1} & \rmif i<n \\
                                abc_1 & \rmif i=n\\ \end{case}
\end{eqnarray*}
where $b$ is an edge from the end of $a$ to the beginning of $c_1$.
The transition matrix (on $\{a,c_1,c_2,\ldots,c_n\}$) is 
$$A=\left( \begin{array}{rrrrcrr}
      1&1&1&0&\ldots&0\\ \hline
      0&0&1&0&\ldots&0\\
      0&0&0&1&\ldots&0\\
      \vdots&\vdots&\vdots&\vdots&\ddots&\vdots\\
      0&0&0&0&\ldots&1\\
      1&1&0&0&\ldots&0\\
 \end{array} \right)$$
The characteristic polynomial of this matrix is $\lambda(\lambda^n-\lambda^{n-1}-2)$, from which we can find the entropy of the system.
In particular $\lambda_{\max}\tendsto1$ as $n\tendsto\infty$, so $\htop\tendsto0$.
\end{example}

%------------------------------------------------------------------------------------------------------------------------------------------------------------

\begin{example}[The effect of boundary components]

Consider the trellis $T_D$ \linebreak shown in \figpartref{xpuncture}{a}.
\fig{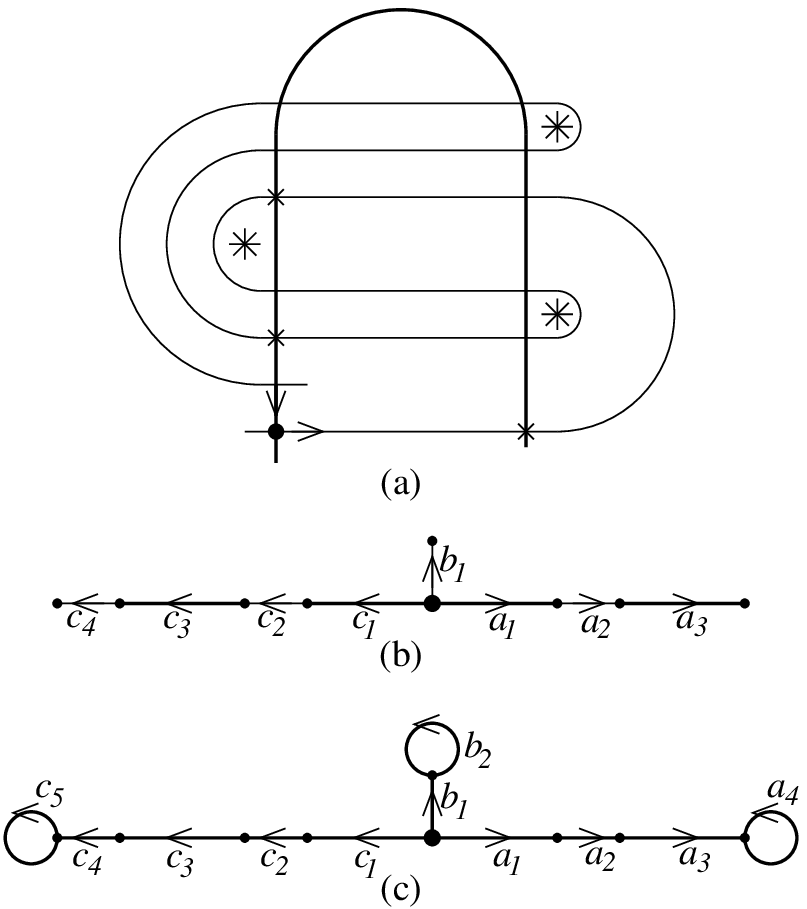}{(a)~Trellis $T_D$ (with puncture points), (b)~Graph $\graph{T_D}$, (c)~Graph $\graph{T_{pA}}$}
The $(U,S,\calO)$-coordinates for the vertices are
\begin{eqnarray*}
(0,0,+),(1,10,-),(2,7,+),(3,4,-),(4,3,+),(5,8,-),(6,9,+), \\
(7,2,-),(8,5,+),(9,6,-),(10,1,+)
\end{eqnarray*}
The graph $\graph{T_D}$ is shown in \figpartref{xpuncture}{c}, and the tight graph map is
$$\begin{array}{llll} a_1\mapsto a_1a_2a_3 & a_2\mapsto\cdot & a_3\mapsto\ba_3\ba_2\ba_1 &  \\
                      b_1\mapsto\cdot & & & \\
                      c_1\mapsto c_1c_2c_3 & c_2\mapsto\cdot & c_3\mapsto\bc_3\bc_2\bc_1 &  c_4\mapsto a_1a_2a_3  \end{array}$$
This map has entropy $\htop=\log2$.

Now suppose the trellis is embedded in a surface with three holes positioned at the stars in \figpartref{xpuncture}{a}
The graph is of the trellis is shown in \figpartref{xpuncture}{c}.
The tight map is  
$$\begin{array}{lllll} a_1\mapsto a_1a_2a_3 & a_2\mapsto a_4 & a_3\mapsto\ba_3\ba_2\ba_1 & a_4\mapsto b_1b_2\bb_1 & \\
                       b_1\mapsto c_1c_2c_3 & b_2\mapsto c_4c_5\bc_4 & & & \\
                       c_1\mapsto c_1c_2c_3 & c_2\mapsto c_4c_5\bc_4 & c_3\mapsto\bc_3\bc_2\bc_1 &  c_4\mapsto a_1a_2a_3 & c_5\mapsto a_4 \end{array}$$
Since the map does not fold of the edge paths $a_1a_2a_3$ and $c_1c_2c_3c_4$, the dynamics of this map are the same as that of $a\mapsto a\ba b$, $b\mapsto c$ and $c\mapsto c\bc a$. 
From this we can show that the characteristic polynomial of the transition matrix has a factor $\lambda^2-3\lambda+1$, from which we obtain entropy $\htop(f)\ge\htop(g_T)=\log(\frac{3+\sqrt{5}}{2})$.

Note that this entropy is larger than that for the trellis in a surface without holes.
Collapsing the holes to points, we obtain a periodic orbit of period $3$. 
The braid type of this orbit is pseudo-Anosov, and the minimal representative has entropy $\log(\frac{3+\sqrt{5}}{2})$, the same as that computed above.
Further, the trellis is exhibited by a blow-up of the pseudo-Anosov homeomorphism.
Thus all the dynamics are forced by the isotopy class in the surface.
\end{example}

%------------------------------------------------------------------------------------------------------------------------------------------------------------

\begin{example}[A toral Anosov trellis]

Let $A$ be that matrix 
$$A=\left( \begin{array}{cc} 2 & 1 \\ 1 & 1 \\ \end{array}\right) $$
The eigenvalues of $A$ are $\half(3\pm\sqrt{5})$ and the eigenvectors are
$$v_u=\left( \begin{array}{cc} 1\\ \frac{-1+\sqrt{5}}{2}\end{array}\right) \hgap v_s\left( \begin{array}{c} -1 \\ \frac{1+\sqrt{5}}{2} \\ \end{array}\right) $$

The trellis $T_A$ of \figpartref{xanosov}{a} occurs in the toral Anosov map with matrix $A$
\fig{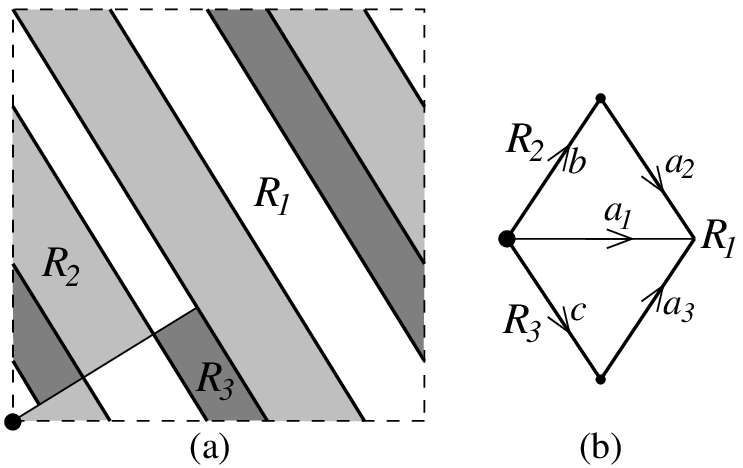}{(a)~Trellis $T_A$, (b)~Graph $\graph{T_A}$}
The points of intersection have coordinates 
\begin{eqnarray*}
q_0 & = & (0,0) \\
q_1 & = & \textstyle{\frac{1}{10}} (-15+7\sqrt{5},25-11\sqrt{5}) \\
q_2 & = & \textstyle{\frac{1}{10}} (- 5+3\sqrt{5},10- 4\sqrt{5}) \\
q_3 & = & \textstyle{\frac{1}{10}} (-10+6\sqrt{5},20- 8\sqrt{5}) \\
q_4 & = & \textstyle{\frac{1}{10}} (    2\sqrt{5}, 5-  \sqrt{5})
\end{eqnarray*}
and the Anosov map $f$ fixes $q_0$ and maps $q_1\mapsto q_2\mapsto q_4$.

The graph $\graph{T_A}$ for $T_A$ is shown in \figpartref{xanosov}{b} and has edges which map:
$$a_1\mapsto a_1,\;a_2\mapsto ba_2,\;a_3\mapsto ca_3,\;b\mapsto ba_2\ba_3\bc \rmand c\mapsto a_1\ba_2\bb$$
If $\alpha=a_1$, $\beta=ba_2$ and $\gamma=ca_3$, then we have
$$\alpha\mapsto\alpha,\;\beta\mapsto\beta\bar{\gamma}\beta \rmand \gamma\mapsto\alpha\bar{\beta}\gamma$$.
Thus the growth rate of the number of periodic points is simply the Perron-Frobenius eigenvalue $\half(3+\sqrt{5})$ of $A$, and all orbits of the Anosov map persist under homotopies preserving the trellis structure.
\end{example}

%------------------------------------------------------------------------------------------------------------------------------------------------------------

\begin{example}[A heteroclinic trellis]

The heteroclinic trellis $T_H$ shown in \figpartref{xheteroclinic}{a} occurs in the Smale horseshoe.
\fig{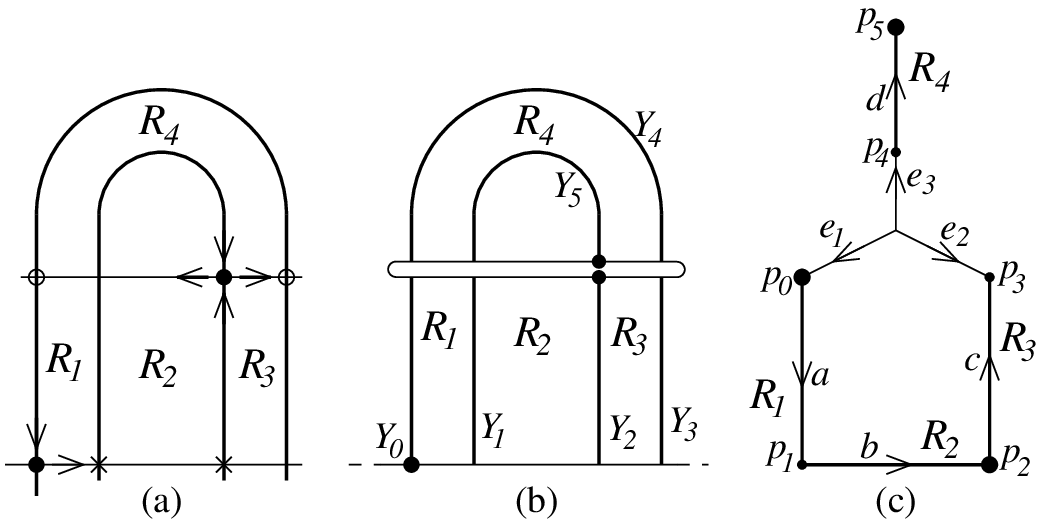}{(a)~Heteroclinic trellis $T_H$, (b)~Surface $\cutt{T_H}$, (c)~Graph $\graph{T_H}$}
There are two saddle fixed points, $p_0$ and $p_1$.
Cutting along the unstable manifold, we obtain the surface $\cutt{T_H}$ of \figpartref{xheteroclinic}{b}, and we can retract this to the graph $\graph{T_H}$ as shown in \figpartref{xheteroclinic}{c}.
The action on the distinguished vertex set is:
$$p_0,p_4\mapsto p_0,\;p_1,p_5\mapsto p_2,\;p_2\mapsto p_5 \rmand p_3\mapsto p_4$$.
The regions $R_1$, $R_2$, $R_3$ and $R_4$ are expanding under the tight map, for which
$$a\mapsto ab, \; b\mapsto c\bar{e}_2e_3d, \; c\mapsto \bar{d} \rmand  e\mapsto ab$$
This gives transition matrix (on $\{a,b,c,d\}$)
$$A=\left( \begin{array}{cccc}
      1&1&0&0\\ \hline
      0&0&1&1\\  
      0&0&0&1\\
      1&1&0&0\\
 \end{array} \right)$$
The characteristic polynomial for $A$ is $\lambda(\lambda^3-\lambda^2-\lambda-1)=0$, and the maximum eigenvalue is $\lambda_{\max}\approx 1.839$.
$\log\lambda_{\max}\approx 0.609$, so $h_\top(f)\ge0.609$ for any map with this trellis action.
Note that this entropy bound is less than that obtained from the horseshoe trellis $T_2$.
\end{example}

%------------------------------------------------------------------------------------------------------------------------------------------------------------

\begin{example}[A trellis with tangential intersections]
\label{example:tangential}

Consider the trellis $T_I$ of \figpartref{xtangency}{a} which occurs in bifurcations from the Smale horseshoe and has tangential intersections.
\fig{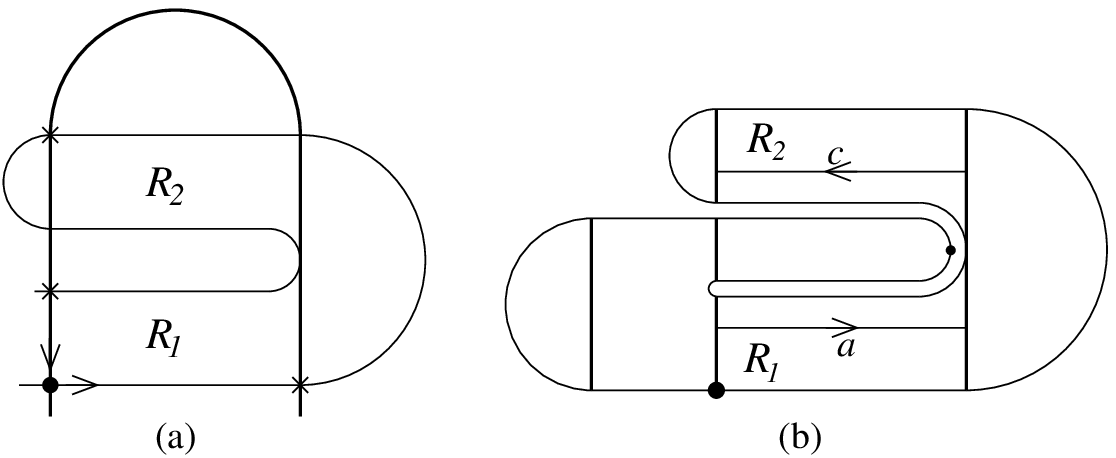}{(a)~Trellis with tangencies $T_I$, (b)~Surface $\cutt{T_I}$}
Cutting along the unstable manifold, we obtain the surface $\cutt{T_I}$ shown in \figpartref{xtangency}{b}
This is not a cross-cut surface, and while there is an exact deformation retract from this surface to a divided graph, we shall study the induced map using the Lefschetz theory.

The cohomology action gives 
$$\alpha\mapsto\alpha+\beta+\gamma, \ \beta\mapsto 0 \rmand \gamma\mapsto -\alpha-\beta-\gamma$$
Just considering the cohomology action on $\alpha$, and $\gamma$, we have Lefschetz matrices
$$A      =\left( \begin{array}{cc} 1& 1\\ \hline -1&-1\\ \end{array} \right) \hgap
  A_{R_1}=\left( \begin{array}{cc} 1& 1\\ \hline  0& 0\\ \end{array} \right) \hgap
  A_{R_2}=\left( \begin{array}{cc} 0& 0\\ \hline -1&-1\\ \end{array} \right) $$
Thus for any word $\regw$ on $R_1$ and $R_2$, $L(A_\regw)=\pm1$ and so $\clperr{\regw}{f}\neq\emptyset$.
Again, we have at least $2^n$ points of period $n$ for $f$, and since $R_1$ and $R_2$ are disjoint, we can again deduce that the topological entropy is at least $\log 2$.
\end{example}

%************************************************************************************************************************************************************

\section{Further Study}
\label{sec:Further Study}

In this paper we describe a general framework for studying maps with tangles.
However there are still many unanswered questions and opportunities for further work.

One particularly important problem is that of optimality of these methods.
This is intimately related to the conditions we place on the map itself.
As an example, consider a homoclinic trellis on the sphere with two intersections, and a map $f$ with this trellis.
If $f$ is a diffeomorphism, we know that $f$ must have a horseshoe in some iterate, and hence be chaotic and have exponential growth of periodic points.
Unfortunately, as previously remarked, we cannot find a lower bound for topological entropy, even though we know if must be strictly positive.
Using the pruning theory of de Carvalho \cite{deCarvalhoPP} we can show that there is a homeomorphism with this trellis with zero entropy.
This homeomorphism has stable and unstable curves at the fixed point, but this fixed point is not hyperbolic.
Therefore, it is not surprising that our methods do not give periodic orbits when applied in this case.

For many examples, we can show that there is a uniformly hyperbolic diffeomorphism with the given trellis which realises the entropy bound given by the asymptotic Nielsen number. As remarked above, this cannot be true in general, but a nice result would be the following

\begin{conjecture}
Let $f$ be a trellis map for the trellis $T$.
Then $N_\infty(f)$ is a lower bound for all maps with trellis $T$ homotopic to $f$.
Further, there is a homeomorphism homotopic to $f$ with topological entropy $N_\infty(f)$, and for all $\epsilon>0$ there is a uniformly hyperbolic diffeomorphism homotopic to $f$ such that $\htop h<N_\infty(f)+\epsilon$.
\end{conjecture}

A possible way of constructing these diffeomorphisms is by using a tight graph map.
For this method to work, we probably need to show that for any trellis map $f$, there is a tight graph map isomorphic to $\cut{f}$ (for a suitable regional decomposition).
Since we cannot in general find a morphism in the category of dynamical systems from a general graph map to a tight one without losing entropy, this could be a tricky problem.

Another interesting problem is the case of non-invertible maps.
We have shown that there are no major problems unless points not in $T^U$ maps over $T^U$, in which case our method breaks down.
Sander \cite{SanderPP} showed that in general, non-invertible maps may have non-trivial tangles but still be non-chaotic.
However, we still may be able to deduce chaos in more general situations than those described here.

Ultimately, we would like to refine this procedure into an algorithm suitable for implementation on a computer.
This requires a way of encoding the important properties of trellises and trellis maps combinatorially.
As we have seen, the $(U,S,\calO)$ coordinate description for the vertices provides a good description of a homoclinic trellis on a sphere; in more complicated cases we have to take into account the homotopy classes of the curves in the surface $M$, and also the way different curves wind round each other.

Having obtained a complete description of a single trellis, we would then like to consider bifurcation sequences.
This requires an especially good understanding of trellises with tangential intersections.
Since Nielsen classes are open in the set of periodic points of a given period, they cannot be removed by sufficiently small perturbations, even if the trellis is destroyed.
Therefore, our analysis of the trellis in \exampleref{tangential} shows that all periodic horseshoe orbits are present at the bifurcation of the trellis, and therefore, given a sufficiently small perturbation, all such orbits of sufficiently low period remain.
However, the possible orderings in which periodic orbits may be destroyed is unknown, though some results have been obtained by Hall \cite{Hall94}.

%************************************************************************************************************************************************************

\providecommand{\bysame}{\leavevmode\hbox to3em{\hrulefill}\thinspace}

%------------------------------------------------------------------------------------------------------------------------------------------------------------


\begin{thebibliography}{Bro71}

\bibitem[BH95]{BestvinaHandell95}
Mladen Bestvina and Michael Handel, \emph{Train-tracks for surface
  homeomorphisms}, Topology \textbf{34} (1995), no.~1, 109--140.

\bibitem[Bro71]{Brown71}
Robert Brown, \emph{The {L}efschetz fixed point theorem}, Scott, Foresman and
  Company, 1971.

\bibitem[Col]{CollinsPPa}
Pieter Collins, \emph{Relative periodic point theory}, Unpublished.

\bibitem[dC]{deCarvalhoPP}
Andr\'e de~Carvalho, \emph{Pruning fronts and the formation of horseshoes},
  Preprint.

\bibitem[Eas86]{Easton86}
Robert Easton, \emph{Trellises formed by stable and unstable manifolds in the
  plane}, Trans. Amer. Math. Soc. \textbf{294} (1986), no.~2, 719--732.

\bibitem[FM93]{FranksMisiurewicz93}
John Franks and Michael Misiurewicz, \emph{Cycles for disk homeomorphisms and
  thick trees}, Nielsen Theory and Dynamical Systems, Contemporary Mathematics,
  1993.

\bibitem[Hal94]{Hall94}
Toby Hall, \emph{The creation of horseshoes}, Nonlinearity \textbf{7} (1994),
  no.~3, 861--924.

\bibitem[KH95]{KatokHasselblatt95}
Anatole Katok and Boris Hasselblatt, \emph{Introduction to the modern theory of
  dynamical systems}, Encyclopedia of Mathematics and its Applications, no.~54,
  Cambridge University Press, 1995.

\bibitem[KO98]{KrauskopfOsinga98}
Bernd Krauskopf and Hinke Osinga Bernd~Krauskopf, \emph{Growing $1$d and quasi-$2$d unstable
  manifolds of maps}, J. Comput. Phys. \textbf{146} (1998), no.~1, 404--419.

\bibitem[Los93]{Los93}
J\'er\^ome~E. Los, \emph{Pseudo-{A}nosov maps and invariant train tracks in the
  disc: A finite algorithm}, Proc. London Math. Soc. (3) \textbf{66} (1993),
  no.~2, 400--430.

\bibitem[San]{SanderPP}
Evelyn Sander, \emph{Homoclinic tangles for noninvertible maps}, Preprint.

\bibitem[Szy95]{Szymczak95FUNDM}
Andrzej Szymczak, \emph{The {C}onley index for decompositions of isolated
  invariant sets}, Fund. Math. \textbf{148} (1995), no.~1, 71--90.

\end{thebibliography}
\end{document}